\documentclass[12pt]{amsart}
\usepackage{amsmath,amssymb,amsbsy,amsfonts,latexsym,amsopn,amstext,
                                               amsxtra,euscript,amscd,bm}
\usepackage{url}
\usepackage[colorlinks,linkcolor=blue,anchorcolor=blue,citecolor=blue]{hyperref}

\usepackage{color}

\usepackage{float}

\hypersetup{breaklinks=true}

\begin{document}

\newtheorem{theorem}{Theorem}
\newtheorem{lemma}[theorem]{Lemma}
\newtheorem{claim}[theorem]{Claim}
\newtheorem{cor}[theorem]{Corollary}
\newtheorem{prop}[theorem]{Proposition}
\newtheorem{definition}{Definition}
\newtheorem{question}[theorem]{Question}
\newcommand{\hh}{{{\mathrm h}}}

\numberwithin{equation}{section}
\numberwithin{theorem}{section}
\numberwithin{table}{section}

\def\sssum{\mathop{\sum\!\sum\!\sum}}
\def\ssum{\mathop{\sum\ldots \sum}}
\def\dsum{\mathop{\sum \sum}}
\def\iint{\mathop{\int\ldots \int}}
\def\li{\mathop{li}}

\def\squareforqed{\hbox{\rlap{$\sqcap$}$\sqcup$}}
\def\qed{\ifmmode\squareforqed\else{\unskip\nobreak\hfil
\penalty50\hskip1em\null\nobreak\hfil\squareforqed
\parfillskip=0pt\finalhyphendemerits=0\endgraf}\fi}

\newfont{\teneufm}{eufm10}
\newfont{\seveneufm}{eufm7}
\newfont{\fiveeufm}{eufm5}
%
%
\newfam\eufmfam
     \textfont\eufmfam=\teneufm
\scriptfont\eufmfam=\seveneufm
     \scriptscriptfont\eufmfam=\fiveeufm
%
%
\def\frak#1{{\fam\eufmfam\relax#1}}

\newcommand{\bflambda}{{\boldsymbol{\lambda}}}
\newcommand{\bfmu}{{\boldsymbol{\mu}}}
\newcommand{\bfxi}{{\boldsymbol{\xi}}}
\newcommand{\bfrho}{{\boldsymbol{\rho}}}

\def\fK{\mathfrak K}
\def\fT{\mathfrak{T}}

\def\fA{{\mathfrak A}}
\def\fB{{\mathfrak B}}
\def\fC{{\mathfrak C}}

\def \balpha{\bm{\alpha}}
\def \bbeta{\bm{\beta}}
\def \bgamma{\bm{\gamma}}
\def \blambda{\bm{\lambda}}
\def \bchi{\bm{\chi}}
\def \bphi{\bm{\varphi}}
\def \bpsi{\bm{\psi}}

\def \wtN {\widetilde N}

\def\vec#1{\mathbf{#1}}


\def\cA{{\mathcal A}}
\def\cB{{\mathcal B}}
\def\cC{{\mathcal C}}
\def\cD{{\mathcal D}}
\def\cE{{\mathcal E}}
\def\cF{{\mathcal F}}
\def\cG{{\mathcal G}}
\def\cH{{\mathcal H}}
\def\cI{{\mathcal I}}
\def\cJ{{\mathcal J}}
\def\cK{{\mathcal K}}
\def\cL{{\mathcal L}}
\def\cM{{\mathcal M}}
\def\cN{{\mathcal N}}
\def\cO{{\mathcal O}}
\def\cP{{\mathcal P}}
\def\cQ{{\mathcal Q}}
\def\cR{{\mathcal R}}
\def\cS{{\mathcal S}}
\def\cT{{\mathcal T}}
\def\cU{{\mathcal U}}
\def\cV{{\mathcal V}}
\def\cW{{\mathcal W}}
\def\cX{{\mathcal X}}
\def\cY{{\mathcal Y}}
\def\cZ{{\mathcal Z}}
\newcommand{\rmod}[1]{\: \mbox{mod} \: #1}

\def\cg{{\mathcal g}}

\def\G{\mathbf G}

\def\e{{\mathbf{\,e}}}
\def\ep{{\mathbf{\,e}}_p}
\def\eq{{\mathbf{\,e}}_q}

\def\em{{\mathbf{\,e}}_m}

\def\Tr{{\mathrm{Tr}}}
\def\Nm{{\mathrm{Nm}}}

 \def\SS{{\mathbf{S}}}

\def\lcm{{\mathrm{lcm}}}

\def\({\left(}
\def\){\right)}
\def\l|{\left|}
\def\r|{\right|}
\def\fl#1{\left\lfloor#1\right\rfloor}
\def\rf#1{\left\lceil#1\right\rceil}
\def\flq#1{\langle #1 \rangle_q}

\def\mand{\qquad \mbox{and} \qquad}

\newcommand{\comml}[1]{\marginpar{%
\begin{color}{blue}
\vskip-\baselineskip 
\raggedright\footnotesize
\itshape\hrule \smallskip #1\par\smallskip\hrule\end{color}}}

\newcommand{\commC}[1]{\marginpar{%
\begin{color}{magenta}
\vskip-\baselineskip 
\raggedright\footnotesize
\itshape\hrule \smallskip C: #1\par\smallskip\hrule\end{color}}}




\hyphenation{re-pub-lished}

\mathsurround=1pt

\def\bfdefault{b}
\overfullrule=5pt

\def \F{{\mathbb F}}
\def \K{{\mathbb K}}
\def \Z{{\mathbb Z}}
\def \Q{{\mathbb Q}}
\def \R{{\mathbb R}}
\def \C{{\\mathbb C}}
\def\Fp{\F_p}
\def \fp{\Fp^*}

\def\Smn{S_{k,\ell,q}(m,n)}

\def\Kmn{\cK_p(m,n)}
\def\psmn{\psi_p(m,n)}

\def\SM{\cS_{k,\ell,q}(\cM)}
\def\SMN{\cS_{k,\ell,q}(\cM,\cN)}
\def\SAMN{\cS_{k,\ell,q}(\cA;\cM,\cN)}
\def\SABMN{\cS_{k,\ell,q}(\cA,\cB;\cM,\cN)}

\def\SIJq{\cS_{k,\ell,q}(\cI,\cJ)}
\def\SAJq{\cS_{k,\ell,q}(\cA;\cJ)}
\def\SABJq{\cS_{k,\ell,q}(\cA, \cB;\cJ)}

\def\sM{\cS_{k,q}^*(\cM)}
\def\sMN{\cS_{k,q}^*(\cM,\cN)}
\def\sAMN{\cS_{k,q}^*(\cA;\cM,\cN)}
\def\sABMN{\cS_{k,q}^*(\cA,\cB;\cM,\cN)}

\def\sIJq{\cS_{k,q}^*(\cI,\cJ)}
\def\sAJq{\cS_{k,q}^*(\cA;\cJ)}
\def\sABJq{\cS_{k,q}^*(\cA, \cB;\cJ)}
\def\sABJp{\cS_{k,p}^*(\cA, \cB;\cJ)}

 \def \xbar{\overline x}

\title[Connected components of power maps]{Connected components of the graph generated by 
power maps  in prime finite fields}

 \author[C. Pomerance]{Carl Pomerance}
\address{Mathematics Department, 
Dartmouth College,
Hanover, NH 03755, USA}
\email{carl.pomerance@dartmouth.edu}

 \author[I. E. Shparlinski] {Igor E. Shparlinski}
\address{Department of Pure Mathematics, University of New South Wales,
Sydney, NSW 2052, Australia}
\email{igor.shparlinski@unsw.edu.au}

\dedicatory{For Jeffrey Outlaw Shallit on his 60th birthday}

\begin{abstract} 
Consider the power pseudorandom-number generator in a finite field $\F_q$.
That is, for some integer $e\ge2$, one considers the sequence $u,u^e,u^{e^2},\dots$
in $\F_q$ for a given seed $u\in \F_q^\times$.  This sequence is eventually periodic.
One can consider the number of cycles that exist as the seed $u$ varies over
$\F_q^\times$.  This is the same as the number of cycles in the  functional graph of the map $x\mapsto x^e$
in $\F_q^\times$.
We prove some estimates for the maximal and average number of
cycles in the case of prime finite fields.
 \end{abstract}
%

\maketitle

\section{Introduction}
\subsection{Set up}

For a prime power $q$, we use $\F_q$ to denote the finite field of $q$ elements.
For a fixed integer $e\ge 2$ we denote by $\cG_{e,q}$ the functional graph of the map 
$x \mapsto x^e$  with vertices formed by the elements of $\F_q^\times$. 
We also denote by   $N(e,q)$ the total number of cycles in  $\cG_{e,q}$.
Alternatively, $N(e,q)$  can be defined as the number of connected components
of $\cG_{e,q}$ when it is considered as an undirected graph. 

By a result of~\cite[Theorem~1]{ChouShp} for prime fields (see also~\cite{VaSha} for $e=2$), which can easily be extended to arbitrary 
finite fields, we have   
\begin{equation}
\label{eq:Neq}
N(e,q) = \sum_{d\mid \rho}\frac{\varphi(d)}{\ell_e(d)}, 
\end{equation}
where $\rho$ is the largest divisor of $q-1$ which is relatively prime to $e$ 
and, for $a,b$ relatively prime and $b$ positive, $\ell_a(b)$  denotes
the multiplicative order of $a$  modulo $b$. 

Here we are interested in the extreme and average values of $N(e,q)$ 
when $e$ is fixed and $q$ varies over primes. 

We remark that under the Generalised Riemann Hypothesis, 
the orders  $\ell_a(b)$ tend to be large (of magnitude $b$ in a logarithmic scale); 
we refer to~\cite{LiPom}. Hence one expects that for most primes
we have $N(e,p) \le p^{o(1)}$. On the other hand, we show that the average 
value of $N(e,p)$ is quite large. 

\subsection{Notation}

Throughout the paper,  the letters  $p$  and $r$  always denote  prime 
numbers while the letters  $a$, $e$, $k$, $m$, and $n$ denote positive integers. 

As usual, for a positive real  number $x$ we use  $\pi(x)$ to denote the number of primes $p\le x$.
Furthermore,  for  integers $a$ and $k \ge 1$ we define
$\pi(x;k,a)$ as the number of primes $p\le x$ in the arithmetic progression
$p \equiv a \pmod k$.  

We also use $P(k)$ and $\varphi(k)$ to denote the largest prime divisor
and the Euler function of $k$, respectively, with $P(1)=1$.

We recall that the    statements $U=O(V)$, $V \gg U$ and $U\ll V$ are all equivalent to the 
inequality $|U|\le cV$ with some positive constant $c$.  In this note, implied constants may depend
 on the exponent  $e$ unless stated otherwise. 
 
\subsection{New results}

First,  we show that $N(e,p)$ is rather large for infinitely many primes $p$. 

\begin{theorem}
\label{thm:Nep Large} For any fixed integer $e\ge 2$,  there are  infinitely many primes $p$ with
$$
N(e,p)\ge p^{5/12 + o(1)}.
$$
\end{theorem}

We also show the following lower bound on the average value of $N(e,p)$.

\begin{theorem}
\label{thm:Nep Aver} For any fixed integer $e\ge 2$ and all sufficiently large real numbers $x$, we have 
$$
\frac {1}{\pi(x)}\sum_{p\le x}  N(e,p)\ge x^{0.293}.
$$
\end{theorem}

\section{Preliminaries}

\subsection{Primes in arithmetic progressions}

We need a version of a result of Alford, Granville and Pomerance~\cite[Theorem~2.1]{AGP} .

\begin{lemma} 
\label{lem:PrimeProg} 
For each fixed $\varepsilon > 0$ and sufficiently large $x$, depending on $\varepsilon$,
there is a finite set $\{m_1,  \ldots, m_t\}$ of integers, where $t$ depends only on  $\varepsilon$, and each $m_i > \log x$, with the following property. If $m\le x^{5/12 - \varepsilon}$, and m is not divisible by any of $m_1,  \ldots, m_t$, then we have uniformly 
over integers $a$ with $\gcd(a,m)=1$,  that
$$
\pi(x;m,a) \gg \frac {1}{\varphi(m)}\pi(x)
$$
where the implied constant depends only on $\varepsilon$.
\end{lemma}
\subsection{Shifted primes with prescribed  smoothness}

We also need the following result, which follows from
 the  work of  Baker and Harman~\cite[Theorem~1]{BH}, which improves the estimate in~\cite{F}.
We recall our convention that $r$ always denotes a prime number 

\begin{lemma} 
\label{lem:ShiftPrime} 
There is an absolute positive constant $\kappa$ with the following property.
Let $u>10$,
$$
v =  \frac{\log u}{\log_2 u} , \quad w = v^{1/0.2961} ,
$$
and let
$$
\cQ  = \left\{r \in [w/(\log w)^\kappa,w]~:~ r-1\mid M_v \right\},
$$
where $M_v$  is  the least common multiple of the integers in $[1,v]$.
Then for $u$ sufficiently large,
we have
$$
\# \cQ  \ge  w/(\log w)^\kappa . 
$$
\end{lemma}

\section{Proofs of main results}

\subsection{Proof of Theorem~\ref{thm:Nep Large}}

We fix some integer $e\ge2$ and a real $\varepsilon >0$. 
For a sufficiently large number $K$ we define $x$ by the equation 
$$
e^K = x^{5/12 - \varepsilon}.
$$
Now let  $m_1,  \ldots, m_t$ be as in Lemma~\ref{lem:PrimeProg}. 

Clearly if $\gcd(m_i, e) > 1$ then $m_i  \nmid e^k-1$.
 For each $i$ with $m_i$ coprime to $e$,  we obviously have 
\begin{equation}
\label{eq:large li}
\ell_e(m_i) \gg \log m_i \gg \log \log x, \qquad i =1, \ldots, m.
\end{equation}
Hence for any integer $h\ge 1$ we have at least 
$$
h- \sum_{i=1}^t \(\frac{h}{\ell_e(m_i)} + 1\) =h  + O(h/\log \log x + 1)
$$
integers $k$ in the interval $[K-h, K]$,  which are 
 not divisible by any of the 
multiplicative orders $\ell_e(m_i)$ for which $\gcd(m_i,e)>1$.  Thus $e^k-1$   is not divisible by any of the 
integers $m_i$, $i =1, \ldots, t$. 
 In particular, we can always find $k \in [K-h_0, K]$, 
where $h_0$ depends only on $\varepsilon$, for which  $m_i  \nmid e^k-1$, $i =1, \ldots, t$.
We fix such an integer $k$ and denote $m = e^k-1$.
Thus by Lemma~\ref{lem:PrimeProg} there exists a prime 
\begin{equation}
\label{eq:small p}
p \ll x = e^{12K/(5 -12 \varepsilon)} \ll  m^{12/(5 -12 \varepsilon)}
\end{equation}
with 
$$
p \equiv 1 \pmod m.
$$
Since $\gcd(m,e)=1$, we have $m\mid\rho$, where $\rho$ is the part of $p-1$ coprime to $e$.
Thus, using $\ell_e(m)=k$, we obtain
$$
N(e,p) = \sum_{d\mid \rho}\frac{\varphi(d)}{\ell_e(d)}\ge\frac{\varphi(m)}{k}\gg\frac{\varphi(m)}{\log m}.
$$

Using the minimal order of the Euler function,  see~\cite[Theorem~328]{HardyWright}, 
we thus have 
$$
N(e,p) \gg \frac{m}{\log m\log\log m},
$$ 
which together with~\eqref{eq:small p} and taking into account that $\varepsilon > 0$ is 
arbitrary, concludes the proof.

\subsection{Proof of Theorem~\ref{thm:Nep Aver}}

We  follow the construction from the proof of~\cite[Theorem~1]{PomShp} which in turn 
is based on some ideas of Erd\H os~\cite{Erd} .

We fix some sufficiently small $\varepsilon >0$ and let $x$ be large.  For 
$$
u = x^{5/12 - \varepsilon}
$$
we consider the  set $\cQ$  and parameters $v$ and $w$  as in Lemma~\ref{lem:ShiftPrime}. 
Furthermore,  let  $m_1,  \ldots, m_t$ be as in Lemma~\ref{lem:PrimeProg}. 
Note that~\eqref{eq:large li} guarantees that for each $i=1, \ldots, t$ with $\gcd(m_i,e)=1$
we have  $\ell_e(m_i) > 1$
and thus we can choose  a prime divisor $r_i$ of $\ell_e(m_i)$ (we do not claim nor require 
these primes to be distinct).
We now remove at most $t$ such primes from the set $\cQ$ and denote 
the remaining set by $\cQ^*$. Thus $\# \cQ^* = \# \cQ + O(1)$.  Note too that
$\cQ^*$ contains no prime dividing $e$.

Put
$$
\nu =\fl{\frac{\log u}{\log w}}
$$
and consider the set $\cS$ of all products of $\nu$ distinct primes from $\cQ^*$.
Clearly
\begin{equation}
\label{eq:range of d}
u \ge w^\nu \ge m  \ge  (w/(\log w)^\kappa)^\nu  = u^{1 + o(1)}
\end{equation}
for every $m \in \cS$.

Furrthermore, using Lemma~\ref{lem:ShiftPrime}, an easy calculation shows that 
\begin{equation}
\label{eq:Card S}
\#\cS =\binom{\# \cQ^*}{\nu} = u^{0.7039+o(1)}.
\end{equation}
For every $m \in \cS$ we have
$$
\ell_e(m) \mid   M_v
$$
and so by the prime number theorem, we obtain that
\begin{equation}
\label{eq:small order}
\ell_e(m)  \le \exp((1+ o(1)) v) = u^{o(1)} = x^{o(1)}.
\end{equation}
Recalling the definition  of $\cQ^*$ we see that for any $m \in \cS$ we have $m_i \nmid m$, 
$i =1, \ldots, t$.
By the choice of $u$ and 
the upper bound in~\eqref{eq:range of d} we see that  by Lemma~\ref{lem:PrimeProg}
we have 
$$
\pi(x;m,1) \gg  \frac{1}{\varphi(m)}\pi(x) = x^{1+o(1)} u^{-1}
$$
for every $m \in \cS$. Thus, using~\eqref{eq:Card S} we obtain
\begin{equation}
\label{eq:total}
\sum_{m \in \cS}  \pi(x;m,1) \ge  x^{1+o(1)} u^{-0.2961}.
\end{equation}

Now, let $\cP$ be the union of all primes $p\le x$ with $m \mid p-1$ for some $m\in\cS$. 
Since, by the classical bound  on the divisor function, each prime $p \in \cP$
can come from at most $x^{o(1)}$ integers $m \in \cS$, we obtain from~\eqref{eq:total} that 
\begin{equation}
\label{eq:Card P}
\#\cP \ge  x^{1+o(1)} u^{-0.2961}. 
\end{equation}

For every $p$ with $m \mid p-1$ for some $m \in \cS$, using~\eqref{eq:small order}
and then~\eqref{eq:range of d}, 
we have
$$
N(e,p) \ge  \frac{\varphi(m)}{\ell_e(m)} = m^{1+o(1)} = u^{1+o(1)}. 
$$
Therefore, using~\eqref{eq:Card P},
$$
\sum_{p\le x} N(e,p) \ge \sum_{p\in \cP} N(e,p) \ge u^{1+o(1)}  \#\cP \ge x^{1+o(1)} u^{0.7039}.
$$
Recalling the choice of $u$ and taking $\varepsilon$ to be sufficiently small, we  conclude  the proof.

\section{Further improvements} 

Hypothetically the exponents in Theorems~\ref{thm:Nep Large} and~\ref{thm:Nep Aver}
may be replaced with any fixed number smaller than 1.  This is true for Theorem~\ref{thm:Nep Large}
on the assumption that we have exponent $1+\varepsilon$ in Linnik's theorem; 
that is,
for each integer $k>k_0(\varepsilon)$ and residue class $a \pmod k$ coprime to $k$, the least prime 
in this residue class 
is smaller than $k^{1+\varepsilon}$.  The proof that $N(e,p)>p^{1-\varepsilon}$ for infinitely many
primes $p$ then follows the same lines as our proof of Theorem~\ref{thm:Nep Large}.

To prove a $1-\varepsilon$ analogue of Theorem~\ref{thm:Nep Aver} we need in addition to
the strong Linnik constant as above, the conjecture that in Lemma~\ref{lem:ShiftPrime} we
may replace the number $0.2961$ with $\varepsilon$.  This conjecture of Erd\H os is
known to follow from the Elliott--Halberstam conjecture.  The proof that the average of $N(e,p)$
for $p\le x$ exceeds $x^{1-\varepsilon}$ is then the same as our proof of Theorem~\ref{thm:Nep Aver}.

The above improvements are probably out of reach.
 However, 
there is a  possible  way to achieve more modest improvements of 
Theorems~\ref{thm:Nep Large} and~\ref{thm:Nep Aver}, which is based 
on a combination of  a recent result of Chang~\cite{Chang} 
with a result of Harman~\cite{Harm2}.  For this approach,  one first 
has to verify that  that the exponent $3/4$  
in~\cite[Equation~(1.2)]{Harm2} can be replaced by any constant  $c < 1$, see also
the remark after~\cite[Theorem~1.2]{Harm1}. Then this result can be combined with
 the bound of Chang~\cite[Theorem~10]{Chang}
on the zero-free region of $L$-functions of characters with smooth moduli, where the 
modulus $m$ is chosen to satisfy two properties
\begin{itemize}
\item $m = e^k-1$ where $k$ is an integer with a small value of $\varphi(k)$, that is, with $\varphi(k) = o(k)$;
\item $m$ is not divisible by a Siegel modulus, which can be  achieved via the same argument 
as that used in the proof of Theorem~\ref{thm:Nep Aver}.
\end{itemize}

Combining these ideas with our approach one is likely to be able replace $5/12$ with  $0.472$
and $0.293$ with $0.332$ in Theorems~\ref{thm:Nep Large} and~\ref{thm:Nep Aver} respectively.
We also note that using the moduli  of the form $m= e^k-1$ with  $\varphi(k) = o(k)$ as in the above, together with the 
version of the Linnik theorem given by Chang~\cite[Corollary~11]{Chang} one can obtain 
an alternative proof of Theorem~\ref{thm:Nep Large}.  However this produces a much sparser sequence 
of primes than in the current proof of Theorem~\ref{thm:Nep Large}.

\section{Further results and directions} 

In~\cite[Theorem~2]{KP} lower bounds are given for the order of $e$ modulo the part of $p-1$
coprime to $e$ that translate to upper bounds for $N(e,p)$.  Indeed, we have for any function
$\varepsilon(p)\downarrow 0$ that $N(e,p)<p^{1/2-\varepsilon(p)}$ for almost all primes $p$
and on the generalized Riemann Hypothesis, $N(e,p)<p^{\varepsilon(p)}$ for almost all $p$.
(These normal-order results are in stark contrast to the above extremal and average-order results.)

One can also consider the average cycle length.  For a positive integer $n$, let $\ell_e^*(n)$ denote
the order of $e$ modulo the prime-to-$e$ part of $n$.  The average cycle length is then
$$
C(e,p) = \frac{1}{p-1} \sum_{d\mid p-1} \varphi(d)\ell^*_e(d).
$$
Note that $\ell_e^*(p-1)=\ell_e(\rho)$, so we have
$$
\frac{\varphi(p-1)}{p-1}\ell_e(\rho)\le C(e,p)\le \ell_e(\rho).
$$
One then sees that results on $\ell_e(\rho)$ immediately translate to results on $C(e,p)$.
So, it follows from~\cite[Theorem~2]{KP} that for any $\varepsilon(p)\downarrow 0$,
we have that for almost all primes $p$, $C(e,p)>p^{1/2+\varepsilon(p)}$.  Further,
the average of $C(e,p)$ for $p\le x$ exceeds $x^{0.592}$ for all sufficiently large values of $x$.
And on the Generalised Riemann Hypothesis, the average exceeds $x^{1-\epsilon}$.  An upper bound for the minimal order
of $C(e,p)$ follows
from the proof of Theorem~\ref{thm:Nep Large}.
In particular, we have $C(e,p)< p^{0.472+o(1)}$ for infinitely many primes $p$.

It would be interesting to generalize the results of this paper to arbitrary finite fields, or
perhaps to consider quantities such as
$$
N(e,p^k), \quad k =1,2, \ldots\ .
$$
For example, we can show that for any fixed choice of $e$ and $p$, for infinitely many $k$ we have
\begin{equation}
\label{eq:Nepk}
N(e,p^k)>\exp(k^{c/\log \log k}),
\end{equation}
where $c$ is a positive constant.  
Indeed, from~\cite[Theorem~1]{EPS}
there are infinitely many positive integers $m$ with $\lambda(m)\le(\log m)^{O(\log \log \log m)}$,
where $\lambda(m)$ is the maximum order of an element in $(\Z/m\Z)^\times$.  Further,
with an easy argument, one can insure that $m$ is coprime to $ep$.  Let
$k=\ell_p(m)\le\lambda(m)$.  We have
$$
N(e,p^k) \ge \varphi(\rho)/\ell_e(\rho) \ge 
\varphi(\rho)/\lambda(\rho) \ge m/\lambda(m),
$$
using~\cite[Lemma~2]{FPS}.  Hence $N(e,p^k)\ge m^{1+o(1)}$.  The small size of $\lambda(m)$ in
comparison to $m$ implies that $m$ is large in comparison to $\lambda(m)$.  In particular, we have
 $m\ge \exp(\lambda(m)^{c/\log \log \lambda(m)})$ for some $c>0$.
 The bound~\eqref{eq:Nepk} follows using $\lambda(m)\ge k$.

It also may be of interest to study the number of cycles of the power generator in the ring
$\Z/n\Z$, where a seed is coprime to $n$.  It is likely that the methods of this paper and
of \cite{KP} should be helpful.

\section*{Acknowledgements}

The authors are very grateful to Glyn Harman for some clarifications concerning the 
possible relaxation of the conditions of~\cite[Theorem~1.2]{Harm2}.  They are also grateful
to the referee for a careful reading.

The first-named author was supported in part by NSF grant number DMS-1440140
at the Mathematical Sciences Research Institute.  He thanks MSRI for their hospitality.

The second-named author thanks 
the Max Planck Institute for Mathematics, Bonn, for the generous 
support and hospitality. He 
 was also supported by ARC Grant DP140100118.


\end{document}